\documentclass{elsart}

\usepackage{amsmath,amssymb}
\usepackage{theorem}

\usepackage{float}
\usepackage{graphicx,here}
\usepackage{psfrag}

\theorembodyfont{\rm}
\newtheorem{problem}{Problem}

\renewcommand{\hat}{\widehat}

\renewcommand{\epsilon}{\varepsilon}

\newcommand{\br}{\bm R}

\newcommand{\bm}[1]{{\boldsymbol #1}}
\renewcommand{\hat}{\widehat}

\begin{document}
\begin{frontmatter}
\title{
Method of Fundamental Solutions with\\
 Optimal Regularization Techniques for \\
the Cauchy Problem of the Laplace Equation\\
with Singular Points}
%\thanksref{NSC}}
\if0
\thanks[NSC]{%
The authors gratefully acknowledge the financial support
of the National Science Council of Taiwan through
the grant No.\ NSC96--2811--E--002--047.}
\fi

\author[NTU]{Takemi Shigeta},\qquad
\ead{shigeta@ntu.edu.tw}
\author[NTU]{D.~L. Young\corauthref{young}}
\corauth[young]{Corresponding author.}
\ead{dlyoung@ntu.edu.tw}
%\thanks{Corresponding author.E-mail: dlyoung@ntu.edu.tw}

\address[NTU]{
Department of Civil Engineering and Hydrotech Research Institute \\
National Taiwan University \\
No.~1, Sec.~4, Roosevelt Road, Taipei 10617, Taiwan}

\begin{abstract}
The purpose of this study is to propose a high-accuracy
and fast numerical method
for the Cauchy problem of the Laplace equation.
Our problem is directly discretized by 
the method of fundamental solutions (MFS).
The Tikhonov regularization method stabilizes a numerical solution
of the problem for given Cauchy data with high noises.
The accuracy of the numerical solution depends on
a regularization parameter of the Tikhonov regularization technique
and some parameters of MFS.
The L-curve determines a suitable regularization
parameter for obtaining an accurate solution.
Numerical experiments show that such a suitable regularization
parameter coincides with the optimal one.
Moreover, a better choice of the parameters of MFS is numerically
observed.
It is noteworthy that
 a problem whose solution has singular points can successfully
be solved.
It is concluded that the numerical method proposed in this paper
is effective for a problem with an irregular domain,
singular points, and the Cauchy data with high noises.
\end{abstract}

\end{frontmatter}

\section{Introduction}
Many kinds of inverse problems have recently been studied
in science and engineering.
The Cauchy problem of an elliptic
partial differential equation is a well known inverse problem.
The Cauchy problem of the Laplace equation
is an important problem which
can be applied to inverse problem of electrocardiography \cite{Franzone}.
Onishi {\it et al.}~\cite{Onishi} proposed an iterative method
 for solving the Cauchy problem of the Laplace equation.
This method reduces the original inverse problem
to an iterative process which alternatively solves
two direct problems.
%As a numerical method of the direct problems,
%the boundary element method (BEM) is applied.
This method, called the adjoint method in 
the papers~\cite{Iijima}, \cite{Shirota},
can solve various inverse problems
by applying many kinds of numerical methods
for solving partial differential equations,
such as the finite difference method (FDM),
the finite element method (FEM),
and the boundary element method (BEM).
The convergence of this method for the Cauchy problem of
the Laplace equation has been obtained~\cite{Shigeta}.

The method of fundamental solutions (MFS)
is effective
for easily and rapidly solving the elliptic well-posed direct problems
in complicated domains.
Mathon and Johnston \cite{Mathon} first showed
numerical results obtained by the MFS.
The papers~\cite{Bogomolny}, \cite{Katsurada}
discuss some mathematical theories on the MFS.
Both of the BEM and the MFS are well known boundary methods,
which discretize original problems
based on the fundamental solutions.
The MFS does not require any treatments for the singularity
of the fundamental solution, while the BEM requires singular integrals.
The MFS is a true meshless method, and can easily be extended
to higher dimensional cases.

Wei {\it et al.}~\cite{Wei} applied
the MFS to the Cauchy problems of elliptic equations.
This method uses the source points distributed outside the domain.
The accuracy of numerical solutions depends on
the location of the source points.
They numerically showed the relation between the accuracy
and the radius of a circle where the source points are distributed.
But, the relation between the accuracy and 
the number of source points has not clearly been given, yet.

Many researchers have solved the Cauchy problem by various methods.
However, to our knowledge, the conventional methods cannot solve
a problem whose solution has singular points
(see \cite{Iijima2} for example).

In this paper, we use the MFS to directly discretize the Cauchy problem
of the Laplace equation.
This problem is an ill-posed problem, where
the solution has no continuous dependence on the 
boundary data.
Namely, a small noise contained in the given Cauchy data
has a possibility to affect sensitively on the accuracy of the solution.
The problem is discretized directly by the MFS
and an ill-conditioned matrix equation is obtained.
A numerical solution of the ill-conditioned equation is
unstable.
The singular value decomposition (SVD) can give
an acceptable solution to such an ill-conditioned matrix equation.
The SVD was successfully applied to the MFS for solving
 a direct problem~\cite{Rama}.
Even though we apply the SVD,
we still cannot obtain an acceptable solution
for the case of the noisy Cauchy data.
We use the Tikhonov regularization to obtain
a stable regularized solution of the ill-conditioned equation.
The regularized solution
depends on a regularization parameter.
Then, we need to determine a suitable regularization parameter
to obtain a better regularized solution.
Hansen~\cite{Hansen} suggested the L-curve as
a method for finding the suitable regularization parameter.
It is known that the suitable parameter is the one
 corresponding to a regularized solution
near the ``corner'' of the L-curve.
We can find the corner of the L-curve as a point with the maximum
curvature~\cite{Hosoda}.

Under the assumption of uniform distribution of the source 
and the collocation points,
we will numerically indicate that
a suitable regularized solution obtained by the L-curve
is optimal in the sense that the error is minimized.
We will respectively show the accuracy and 
the optimal regularization parameter
against a noise level.
We will also mention influence of 
%the radius of a circle where the source points are distributed
%and
 the total numbers of the 
source and the collocation points on accuracy.
We will show that our method is effective
for a problem  whose solution has singular points.
It is noteworthy that such kind of problems
can also successfully be solved.

Section~\ref{sec:ps} introduces the Cauchy problem.
In Section~\ref{sec:mfs}, the MFS discretizes the problem.
In Section~\ref{sec:rs}, the singular value decomposition,
the Tikhonov regularization and the L-curve are used
to obtain a suitable regularized solution.
In Section~\ref{sec:ne}, numerical experiments
confirm that the suitable regularization parameter by the L-curve
coincides with the optimal one that minimizes
the error between the regularized solution and the exact one.
The error and the optimal regularization parameter
against the noise level of the Cauchy data are respectively shown.
Then, our interest is how to choose the following three parameters
used in MFS:
the numbers of collocation points,
the number of source points,
and the radius of a circle where source points are
distributed.
A better choice of the parameters is also observed.
A problem with an irregular domain and
a problem whose solution has singular points
are successfully solved, respectively.
Section~\ref{sec:cnc} concludes the paper.

\section{Problem Setting}\label{sec:ps}
We consider the Laplace equation $-\Delta u=0$ in
a two-dimensional bounded domain $\Omega$
enclosed by the boundary $\Gamma$.
We prescribe Dirichlet and Neumann boundary conditions
simultaneously on a part of the boundary $\Gamma$,
denoted by $\Gamma_1$, as follows:
\[ u=f,\quad \frac{\partial u}{\partial n}=g \quad
\text{on}\quad\Gamma_1, \]
where $f$ and $g$ denote given continuous functions
 defined on $\Gamma_1$,
and $n$ the unit outward normal to $\Gamma_1$.
Then, we need to find the boundary value $u$
on the rest of the boundary $\Gamma_2:=\Gamma\setminus\Gamma_1$
or the potential $u$ in the domain $\Omega$.
This problem is called the Cauchy problem of the Laplace equation,
and the boundary data are called the Cauchy data.

Our Cauchy problem is described as follows:
\begin{problem}
For the given Cauchy data $f,g\in C(\Gamma_1)$,
find $u\in C(\Gamma_2)$ or $u\in C^2(\Omega)\cap C^1(\overline\Omega)$
 such that
\begin{alignat}{7}
 -\Delta u&=0 &\quad &\text{in} &\quad &\Omega, \label{eq:cauchylap}\\
u=f, \quad \frac{\partial u}{\partial n}
&=g & \quad & \text{on} & \quad & \Gamma_1.
\label{eq:cauchydata1}
\end{alignat}
\end{problem}

The Cauchy problem is a well known ill-posed problem.
We can show the instability of the solution to
the Cauchy problem of the Laplace equation as follows:
For example,
in the case where
\begin{align*}
\Omega&=(0,1)^2=\{(x,y):0<x<1,\ 0<y<1\},\\
%,\quad 
\Gamma&=[0,1]\times\{0\}=\{(x,0):0\leq x\leq 1\},\\
f(x,0)&=\frac 1{n^k}\sin(nx),\quad
g(x,0)=0 \quad (k>0),
\end{align*}
the solution is given by
\[ u(x,y)
=\frac 1{n^k}\sin(nx)\cosh(ny). \]
Here, we can see that 
\[ \sup_{\bm x\in\Gamma}|f(\bm x)|\to 0,\qquad 
\sup_{\bm x\in\overline\Omega}|u(\bm x)|\to \infty, \]
from which we know that
the solution $u$ of the Cauchy problem
does not depend continuously on
 the Cauchy data $f$ and $g$.

\section{Discretization by the Method of Fundamental Solutions}
\label{sec:mfs}
\begin{figure}[H]
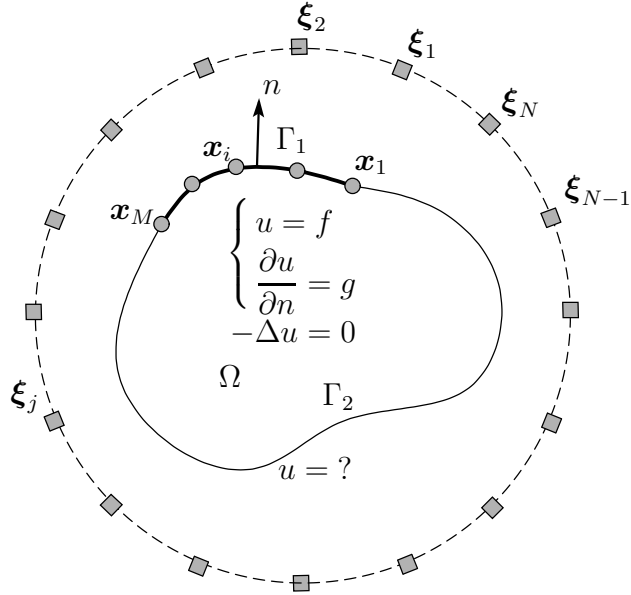

\hspace*{-10ex}
\centering
\input dom31.tex
\caption{Problem setting and distributions of 
collocation and source points}\label{fig:pdp}
\end{figure}
The fundamental solution of the Laplace equation
in two dimensions is defined as
\[ \varphi^*(r):=-\frac 1{2\pi}\ln r \]
for $r=|\bm x|=\sqrt{x^2+y^2}$,
which is a solution to
\[ -\Delta \varphi^*=\delta(\bm x). \]

We distribute the collocation points
$\{\bm x_i\}_{i=1}^M\subset\Gamma_1$ on the boundary
where the Cauchy data are prescribed,
and the source points
$\{\bm \xi_j\}_{j=1}^{N}\subset\overline{\Omega}^c$
along a circle outside the domain (Figure~\ref{fig:pdp}).
We approximate $u$ by $u_{N}$:
\begin{equation} u(\bm x)\approx u_{N}(\bm x)
:=\sum_{j=1}^{N}w_j\varphi_j(\bm x), \label{eq:aprx21}
\end{equation}
where the basis function is defined as
\begin{equation}
 \varphi_j(\bm x):=\varphi^*(|\bm x-\bm \xi_j|) 
\label{eq:basis}
\end{equation}
and $\{w_j\}_{j=1}^N$ are expansion coefficients
to be determined below.
Since the basis functions (\ref{eq:basis}) have no singular points
in $\Omega$,
the approximate function $u_N$ satisfies 
the Laplace equation (\ref{eq:cauchylap}).
Substituting (\ref{eq:aprx21}) into (\ref{eq:cauchydata1})
and assuming that (\ref{eq:cauchydata1}) is satisfied
at the collocation points,
we have
\[ \left\{
\begin{split}
\sum_{j=1}^{N}w_j\varphi_j(\bm x_i)&=f(\bm x_i), \\
\sum_{j=1}^{N}w_j\frac{\partial \varphi_j}{\partial n}
(\bm x_i)&=g(\bm x_i), 
\end{split}
\right. \quad i=1,2,\ldots, M \]
or in the matrix form:
\begin{equation}
A\bm w=\bm b, \label{eq:acb}
\end{equation}
where the matrix $A=(a_{ij})\in\bm R^{2M\times N}$
and the vectors
$\bm w=(w_j)\in\bm R^{N}$, $\bm b=(b_i)\in\bm R^{2M}$
are defined by
\begin{align*}
a_{ij}&:=\left\{ \begin{array}{cl}
\varphi_j(\bm x_i), & i=1,2,\ldots, M \\
\displaystyle
\frac{\partial \varphi_j}{\partial n}(\bm x_{i-M}),
 & i=M+1,M+2,\ldots, 2M
\end{array}\right.,\quad j=1,2,\ldots,N,\\
b_i&:=\left\{ \begin{array}{cl}
f(\bm x_i), & i=1,2,\ldots, M \\
g(\bm x_{i-M}), & i=M+1,M+2,\ldots, 2M
\end{array}\right..
\end{align*}

\section{Regularized Solution}\label{sec:rs}
\subsection{Singular value decomposition}
In general, (\ref{eq:acb}) has the case where the exact solution
$\bm w$
 does not exist in the conventional sense.
As an exact solution,
we consider the least square and least norm solution
$\bm w_0$ defined by 
\begin{equation}
 \|\bm w_0\|=\min_{\bm w\in W}\|\bm w\|,\quad
W:=\{\hat{\bm w}:\|A\hat{\bm w}-\bm b\|=
\min_{\bm w\in\br^N}\|A\bm w-\bm b\|\} .
\label{eq:lslns}
\end{equation}
In this paper, we refer the solution $\bm w_0$
as the exact solution to (\ref{eq:acb}).

For the matrix $A\in\br^{2M\times N}$,
the singular value decomposition (SVD) can be written as follows:
\[ A=U\Sigma V^T=\sum_{i=1}^r\sigma_i\bm u_i\bm v_i^T, \]
where
\begin{gather*}
 U=(\bm u_1\ \bm u_2\ \cdots\ \bm u_{2M})\in\br^{2M\times 2M},\qquad
V=(\bm v_1\ \bm v_2\ \cdots\ \bm v_N)\in\br^{N\times N},\\
\Sigma=(s_{ij})\in\br^{2M\times N},
\qquad s_{ij}=\left\{ \begin{array}{cc}
\sigma_i & (i=j) \\
0 & (i\neq j)
\end{array} \right.,\\
 UU^T=I_{2M}\in\br^{2M\times 2M},\qquad VV^T=I_N\in\br^{N\times N}, \\
 2M\geq N,\qquad \sigma_1\geq\cdots\sigma_r>0,\quad
\sigma_{r+1}=\cdots=\sigma_N=0 
\end{gather*}
with the identity matrices $I_{2M}$ and $I_N$.
The non-negative values $\{\sigma_i\}_{i=1}^N$ are called
the singular values of the matrix $A$.
Using the SVD, we can express the solution to (\ref{eq:lslns}) as
\[ \bm w_0=\sum_{i=1}^r\frac{(\bm u_i,\bm b)}{\sigma_i}\bm v_i. \]

In a real problem, the Cauchy data $f$ and $g$ contain some noises.
We consider the following equation
instead of (\ref{eq:acb}):
\begin{equation}
 A\bm w=\bm b^\delta,\qquad
\bm b^\delta=\bm b+\Delta\bm b
 \label{eq:awbd}
\end{equation}
with
the noise vector $\Delta\bm b\in\bm R^{2M}$.
Then, the solution $\bm w_0^\delta$ to
\begin{equation}
 \|\bm w_0^\delta\|=\min_{\bm w\in W^\delta}\|\bm w\|,\quad
W^\delta:=\{\hat{\bm w}:\|A\hat{\bm w}-\bm b^\delta\|=
\min_{\bm w\in\br^N}\|A\bm w-\bm b^\delta\|\} 
\label{eq:lslns2}
\end{equation}
is quite different from the exact solution $\bm w_0$ to (\ref{eq:lslns})
since the solution is discontinuous for the Cauchy data.
We need to find a good approximation to $\bm w_0$.

\subsection{Tikhonov regularization}
In order to obtain a good approximate solution to (\ref{eq:lslns2}),
we consider minimizing the following functional
with a regularization parameter $\alpha>0$
according to the Tikhonov regularization:
\begin{equation}
 J_\alpha^\delta(\bm w):=\|A\bm w-\bm b^\delta\|^2
+\alpha^2\|\bm w\|^2. \label{eq:fnctnl}
\end{equation}
It is easy to see that the functional $J_\alpha^\delta$
is strictly convex for any $\alpha>0$.
Hence, $J_\alpha^\delta$ has a unique minimum point
 $\bm w_\alpha^\delta$ called the regularized solution:
\[ J_\alpha^\delta(\bm w_\alpha^\delta)=\min_{\bm w\in\br^N}
J_\alpha^\delta(\bm w). \]
We know that $\bm w_\alpha^\delta$ is the solution to
\begin{equation}
 (A^TA+\alpha^2 I_N)\bm w_\alpha^\delta
=A^T\bm b^\delta. \label{eq:fncmateq}
\end{equation}
The equation (\ref{eq:fncmateq}) is uniquely solvable
since the matrix $(A^TA+\alpha^2 I_N)$ is symmetric positive definite.
The SVD of $\bm w_\alpha^\delta$ can be expressed in the form:
\[ \bm w_\alpha^\delta=\sum_{i=1}^r
\gamma_i\frac{(\bm u_i,\bm b^\delta)}{\sigma_i}\bm v_i \]
with the filter factor $\gamma_i:={\sigma_i^2}/({\sigma_i^2+\alpha^2})$.
Then, substituting $\bm w_\alpha^\delta$ into (\ref{eq:aprx21}),
we find the approximate potential $u_N$ in $\Omega\cup\Gamma_2$.

The error between the regularized solution for the noisy data
and the exact solution is decomposed into
\begin{equation}
 \bm w_\alpha^\delta-\bm w_0
=(\bm w_\alpha^\delta-\bm w_\alpha^0)
+(\bm w_\alpha^0-\bm w_0)
=\sum_{i=1}^r\gamma_i
\frac{(\bm u_i,\Delta\bm b)}{\sigma_i}\bm v_i
+\sum_{i=1}^r(\gamma_i-1)\frac{(\bm u_i,\bm b)}{\sigma_i}\bm v_i. 
\label{eq:werr}
\end{equation}
The first term is the perturbation error due to
the relative noise $\Delta\bm b$
and
the second term is the regularization error
caused by regularization of the exact $\bm b$.
When $0<\alpha\ll 1$, we see that $\gamma_i\approx 1$
for most of $i$,
and the error $\bm w_\alpha^\delta-\bm w_0$
is dominated by the perturbation error.
On the other hand, when $\alpha\gg 1$,
we see that $\gamma_i\ll 1$
and the error $\bm w_\alpha^\delta-\bm w_0$
is dominated by the regularization error.
In the next subsection, we will consider a useful method for
 finding a suitable regularization parameter
to minimize both of the perturbation and the regularization errors.

\subsection{L-curve}
To find 
a suitable
regularization parameter,
Hansen~\cite{Hansen} suggests the L-curve,
which is defined as
the continuous curve consisting of all the point
$(\|A\bm w_\alpha^\delta-\bm b^\delta\|,\|\bm w_\alpha^\delta\|)$
for $\alpha>0$:
\[ {\cal L}:=\{(\|A\bm w_\alpha^\delta-\bm b^\delta\|,
\|\bm w_\alpha^\delta\|): \alpha>0 \}. \]
For fixed $\alpha>0$, we get $\bm w_\alpha^\delta$
and then can calculate
the residual norm $\|A\bm w_\alpha^\delta-\bm b^\delta\|$
and the solution norm $\|\bm w_\alpha^\delta\|$.
Thus, the L-curve can be plotted as a set of all the points of
the residual norms as abscissa and the solution norms as ordinate
for all $\alpha>0$.

The L-curve is plotted in double logarithm,
and displays the compromise between
minimization of the perturbation error and the regularization error
in (\ref{eq:werr}).
A suitable regularization parameter is given by
the one corresponding to a regularized solution
near the ``corner'' of the L-curve.
The ``corner'' can be regarded as the point
where the curvature of the L-curve becomes maximum
 \cite{Hansen3}, \cite{Hosoda}.

In the case when $2M=N$ and the Cauchy data $f, g$
have no noises,
if $A$ is non-singular,
we can directly solve (\ref{eq:acb})
to obtain a solution with high accuracy.
However, we cannot guarantee that (\ref{eq:acb}) is always solvable.
Even if there exists the inverse matrix,
the solution to (\ref{eq:acb}) for the noisy Cauchy data differs from
the exact solution.
On the other hand, the regularized solution
 by the Tikhonov regularization is always uniquely determined
for $\alpha>0$.
In the next section,
our numerical experiments will show that
the suitable regularization parameter given by L-curve
coincides with
the optimal one
 $\alpha_{\rm opt}$
defined by
\[ \|\bm w_0-\bm w_{\alpha_{\rm opt}}^\delta\|
=\min_{\alpha>0}\|\bm w_0-\bm w_\alpha^\delta\|. \]

\section{Numerical Experiments}\label{sec:ne}
\subsection{Circular domain}\label{sec:ss}
We first consider a harmonic function $u(x,y)=e^x\cos y-e^y\sin x$
in a unit disk
$\Omega:=\{(x,y):x^2+y^2<1\}$.
According to the exact potential $u$, the exact Cauchy data
are given by
$f=u$ and $g=\partial u/\partial n$
on the fourth part of the whole boundary $\Gamma$,
which is defined by
\[ \Gamma_1:=\{(x,y): x^2+y^2=1,\ x>0,\ y>0\}. \]
We now assume that the exact potential $u$ is unknown,
and identify a boundary value on the rest of the boundary
$\Gamma_2:=\partial\Omega\setminus\Gamma_1$
from the noisy Cauchy data
$f^\delta=(1+\epsilon)f$ and
$g^\delta=(1+\epsilon)g$,
where $\epsilon=\epsilon(x,y)$ is a uniform random noise 
such that
$-\delta\leq \epsilon(x,y) \leq \delta$
with the relative noise level of $100\delta\%$.

We distribute uniformly the collocation points
$\{\bm x_i\}_{i=1}^M\subset\Gamma_1$
and the source points
$\{\bm \xi_j\}_{j=1}^{N}\subset\overline{\Omega}^c$
as follows:
\begin{equation}
 \left\{
\begin{split}
\bm x_i&=(\cos\theta_i,\sin\theta_i),
\quad \theta_i=\frac{2\pi(i-1)}{4M}+\frac{\pi}{4M},
\quad i=1,2,\ldots,M,\\
\bm \xi_j&=(R\cos\hat\theta_j,R\sin\hat\theta_j),
\quad \hat\theta_j=\frac{2\pi(j-1)}{N}+\frac{\pi}{N},
\quad j=1,2,\ldots,N,
\end{split}\right. \label{eq:dispnt}
\end{equation}
where $R>1$ is the radius of the circle where
the source points are distributed.
We adopt the MATLAB code for solving discrete ill-posed problems
based on SVD, made by Hansen~\cite{Hansen2}, \cite{Hansen3}, to
our numerical computations.
Due to the maximum principle,
it is sufficient
to confirm the boundary error between the identified potential $u_N$
and the exact one $u$ rather than the domain error
in our numerical experiments.
We define the maximum relative error on the boundary 
by
\[ e:=\frac{\|u_N-u\|_\infty}{\|u\|_\infty}, \]
where the maximum norm on the boundary denotes
\[ \|u\|_\infty=\sup_{\bm x\in\Gamma}|u(\bm x)|,
\qquad \forall u\in C(\Gamma). \]

\begin{figure}[H]
\centering
\includegraphics[width=8cm]{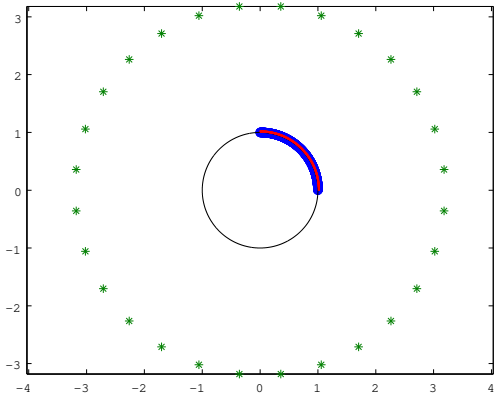}
\caption{Distributions of collocation and source points
($(R,M,N)=(3.2,600,28)$)}
\label{fig:distcs}
\end{figure}
\begin{figure}[H]
\centering
\includegraphics{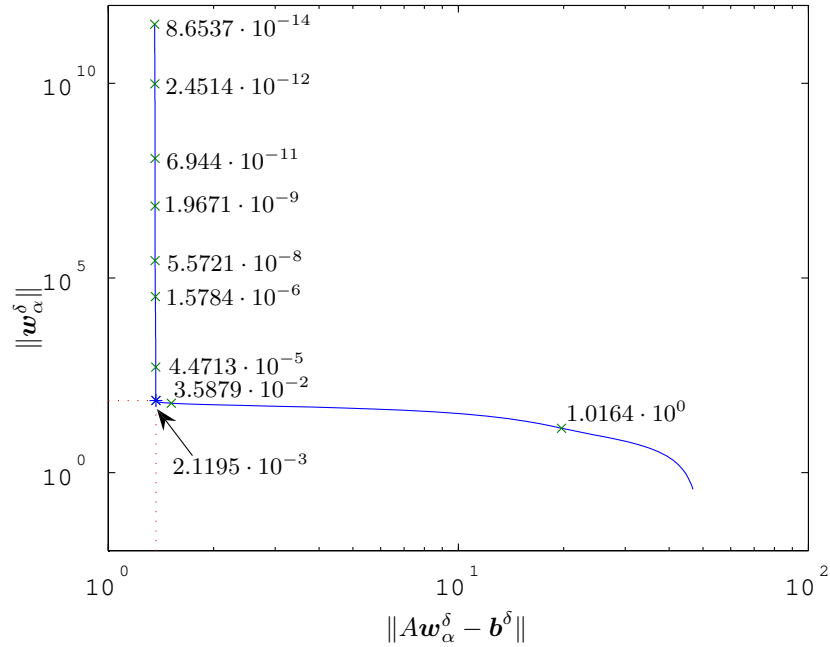}
\caption{L-curve (The corner is located at 
the point for 
$\alpha=2.1195\times 10^{-3}$)}
\label{fig:lcurve}
\end{figure}

\begin{figure}[H]
\centering
%\hspace*{-9em}
\includegraphics{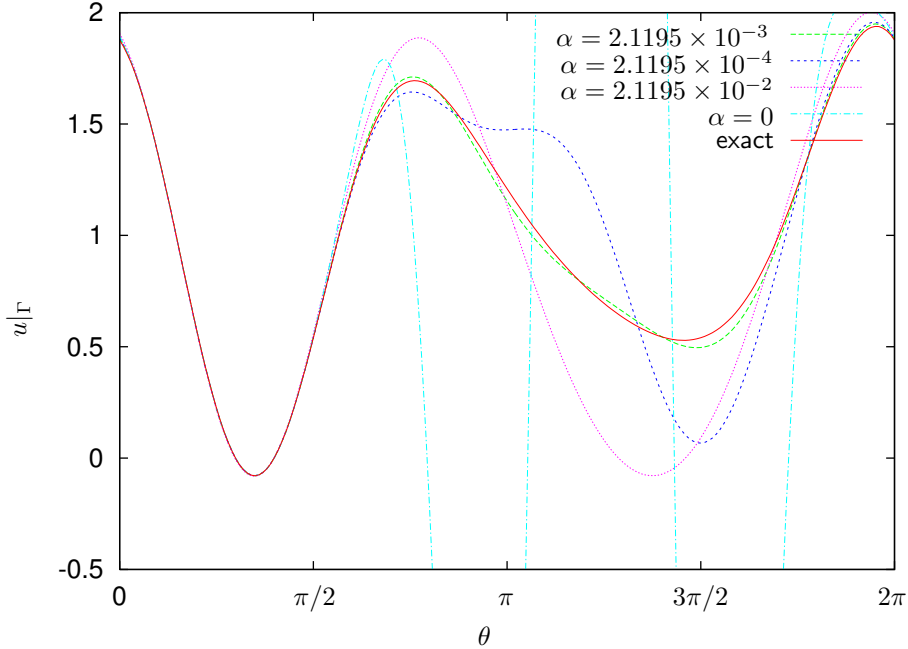}
\caption{The regularized solutions versus the exact solution
($5\%$ noise level)}
\label{fig:regsol}
\end{figure}

\begin{table}[H]
\caption{The maximum relative errors for each $\alpha$}
\label{tab:err}\quad

\centering
\begin{tabular}{c|ccc}
$\alpha$ & $2.12\times 10^{-4}$ & $2.12\times 10^{-3}$
 & $2.12\times 10^{-2}$ \\ \hline
error & 0.2743 & 0.0305 & 0.3382
\end{tabular}
\end{table}

In the first experiment, the relative noise level of the Cauchy data
is assumed to be $5\%$ ($\delta=0.05$).
We set the parameters $(R,M,N)=(3.2,600,28)$.
Figure~\ref{fig:distcs} shows 
the distributions of the collocation and the source points.
As we can see in Figure~\ref{fig:lcurve},
the corner of the L-curve is located at
the point $(\|A\bm w_\alpha^\delta-\bm b^\delta\|,
\|\bm w_\alpha^\delta\|)$
with the regularization parameter $\alpha=2.1195\times 10^{-3}$.
Figure~\ref{fig:regsol} shows the regularized solutions
on the boundary for 
$\alpha=2.1195\times 10^{-4},
2.1195\times 10^{-3},2.1195\times 10^{-2},0$.
We can see that the solution is quite unstable if 
$\alpha=0$, that is, if
the regularization is unapplied.
Comparing the other solutions for $\alpha=2.1195\times 10^{-4},
\alpha=2.1195\times 10^{-3},2.1195\times 10^{-2}$,
we can confirm that $\alpha=2.1195\times 10^{-3}$ is 
a suitable regularization parameter
to obtain a better approximate solution (Table~\ref{tab:err}).

\begin{figure}[H]
\centering
\hspace*{-9em}
\includegraphics{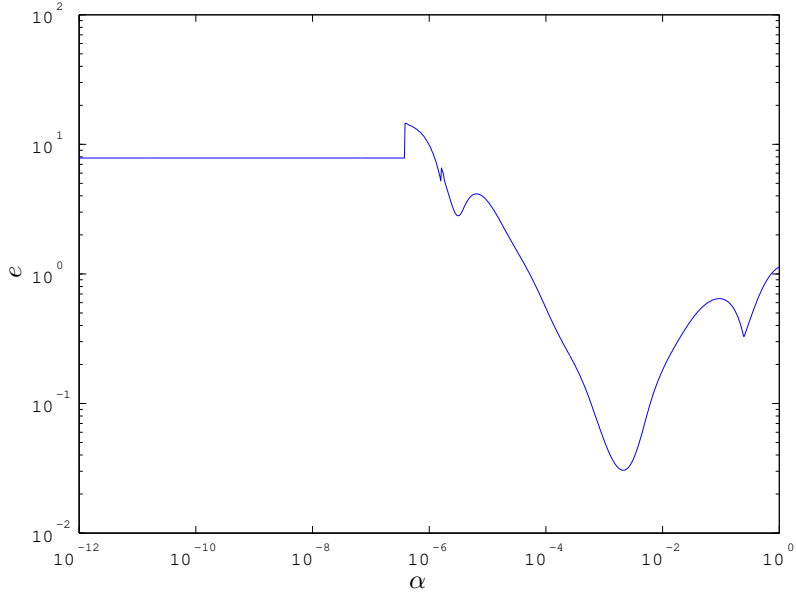}
\caption{Maximum relative error against regularization parameter
}\label{fig:alerr}
\end{figure}
From Figure~\ref{fig:alerr}, we can see that the 
maximum relative error
reaches a minimum at $\alpha=10^{-2.673}\approx 2.12 \times 10^{-3}$,
which coincides with the suitable regularization parameter
obtained by the L-curve.
Hence, we know that
the optimal regularization parameter can be given as the one
corresponding to a regularized solution at the corner of the L-curve.

\begin{figure}[H]
\centering
\includegraphics[width=0.475\textwidth]{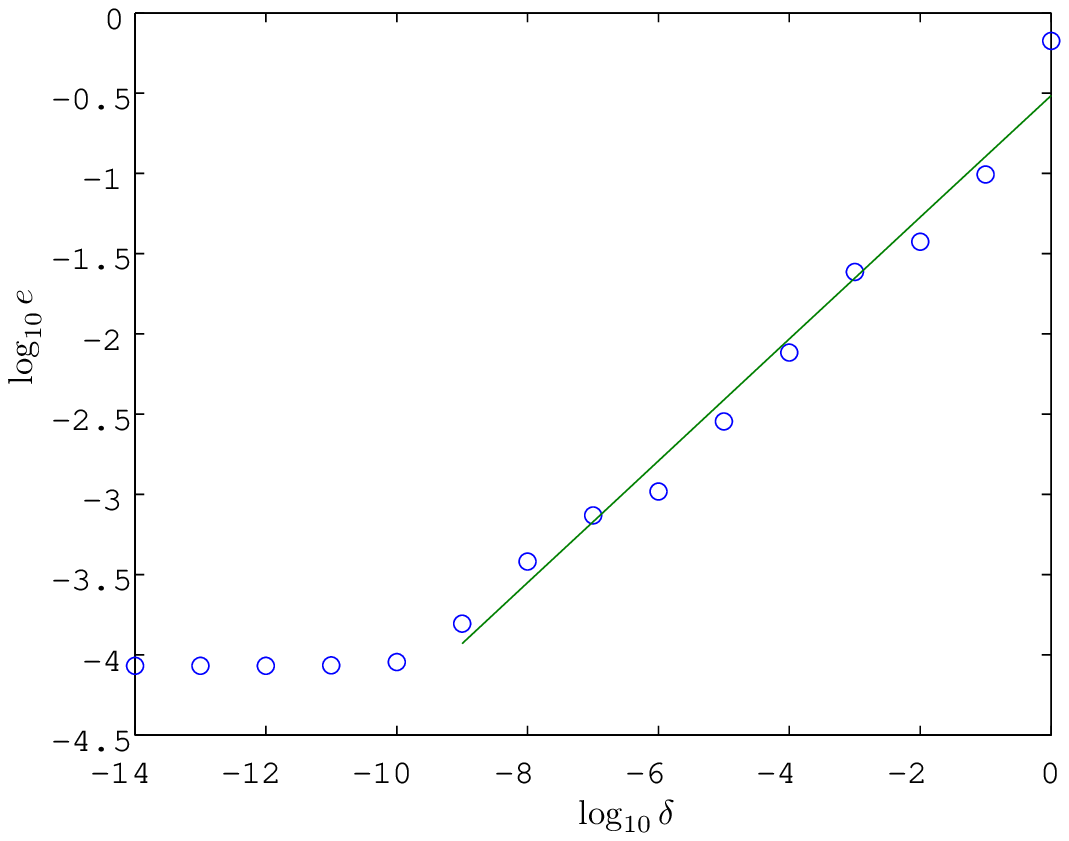}
\hfil
\includegraphics[width=0.475\textwidth]{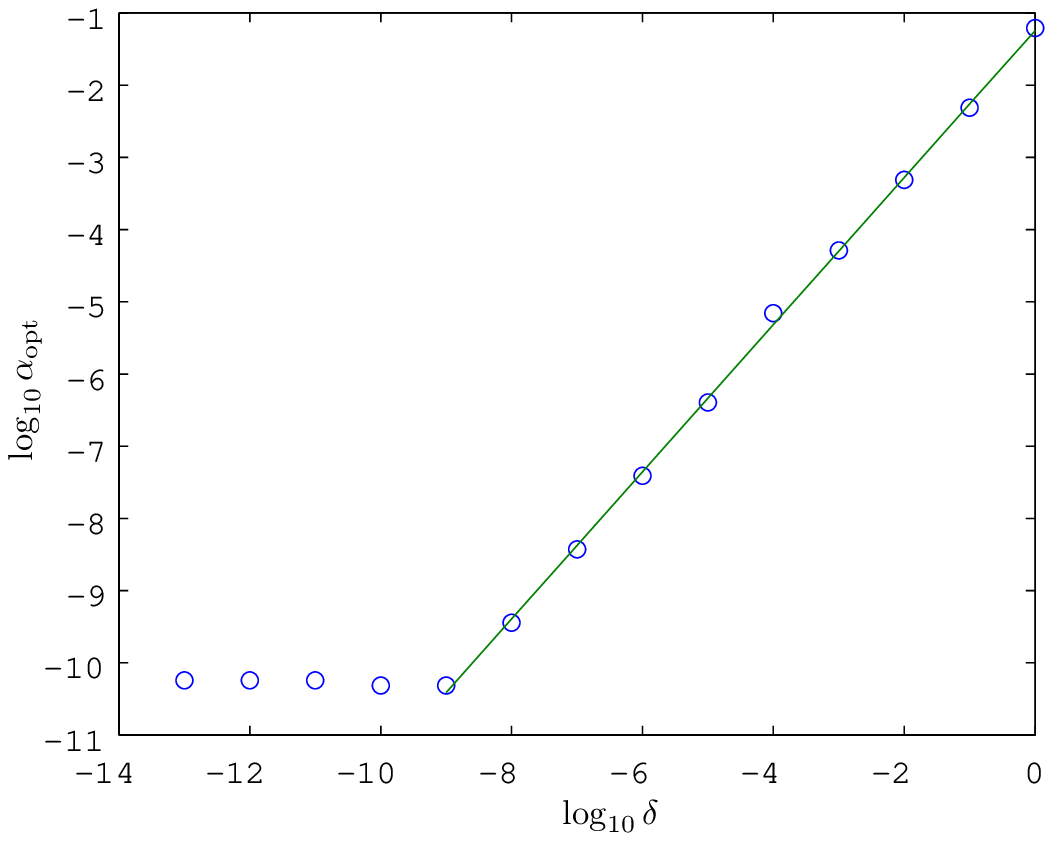}

\centering
(a) \hfil\hfil (b)
\caption{(a) Maximum relative error
 and (b) the optimal regularization parameter
against relative noise level
}\label{fig:er_err_or}
\end{figure}
Figure \ref{fig:er_err_or} (a) shows
the maximum relative error $e$ for the optimal regularized potential
against the relative noise level $\delta$.
The regression line 
in the interval $[-9,0]$ is expressed by
$\log_{10}e=0.37951\log_{10}\delta-0.22672$. 
For the optimal regularized potential $u_N$,
we have 
$e=O(\delta^{0.38})$ for $\delta\geq 10^{-9}$.
Figure \ref{fig:er_err_or} (b) indicates
the optimal regularization parameter $\alpha_{\rm opt}$ against
the relative noise level $\delta$
and
 the regression line 
in the interval $[-9,0]$ given by
$\log_{10}\alpha_{\rm opt}=1.0186\log_{10}\delta-1.2434$.
From this numerical result,
 we can obtain the relation 
$\alpha_{\rm opt}=
O(\delta)$ for $10^{-9}\leq\delta\leq 1$.

After setting the parameters $(R,M,N)$, we can obtain a suitable
regularized solution based on the Tikhonov regularization
and the L-curve.
Now, our problem is how to choose suitable parameters
$(R,M,N)$.

\begin{figure}[H]
\centering
\includegraphics[width=8cm]{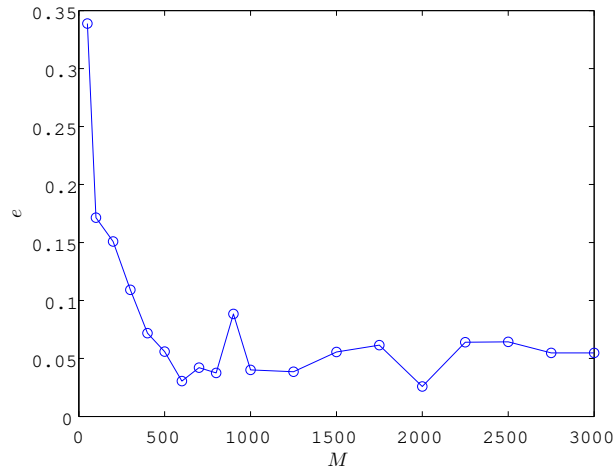}
\caption{Maximum relative error against number of collocation points
}\label{fig:errnc}
\end{figure}

Figure~\ref{fig:errnc} shows the maximum relative error
against the number of collocation points.
We know from this result that we need to take 
sufficiently many collocation points to obtain accurate solutions.

\begin{figure}[H]
\centering
%\includegraphics[width=5cm,height=7cm]{rs_ns_err_global2.eps}
%\includegraphics[width=5cm,height=7cm]{rs_ns_err_local2.eps}
\if0
\includegraphics[width=7.5cm,height=6.5cm]{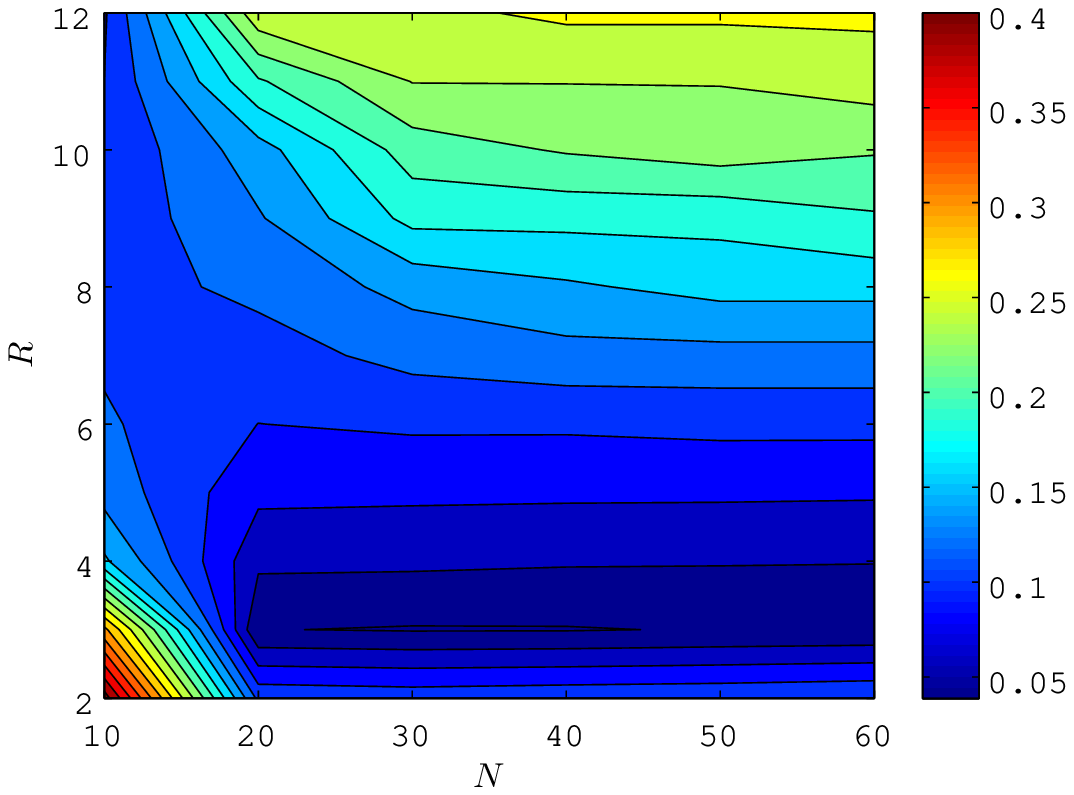}
\hfil
\includegraphics[width=7.5cm,height=6.5cm]{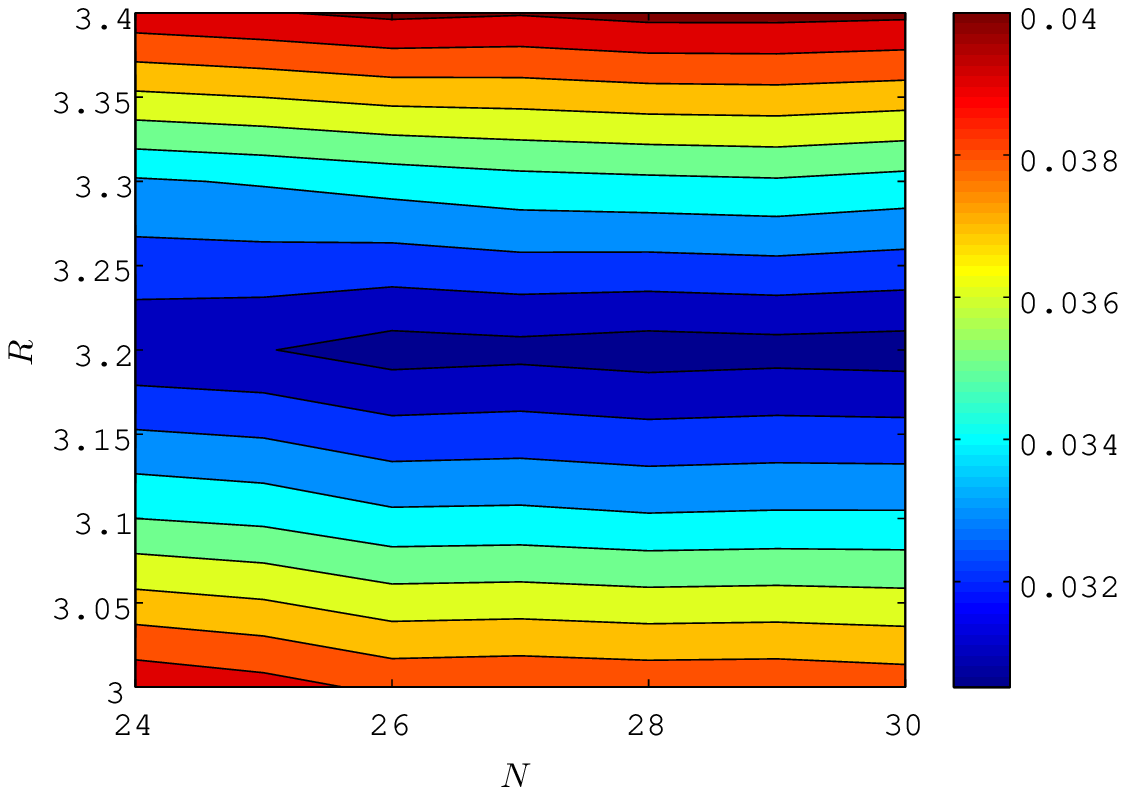}
\fi
\includegraphics[width=0.475\textwidth,height=0.41\textwidth]
{rs_ns_err_global3c.eps}
\hfil
\includegraphics[width=0.475\textwidth,height=0.41\textwidth]
{rs_ns_err_local3c.eps}

(a) \hfil\hfil (b)
\caption{Contour line of the maximum relative error against $(N,R)$
($M=600$);
 (a) $(N,R)\in[10,60]\times[2,12]$,
(b) $(N,R)\in[24,30]\times[3,3.4]$
}\label{fig:mbeRN}
\end{figure}
Figure~\ref{fig:mbeRN} shows the contour line of
the maximum relative error $e$ against $(N,R)$
for the fixed number of collocation points $M=600$.
Through this result, we know that
the maximum relative error is roughly independent of
the number of source points $N$
for the fixed radius $R$ of the circle where source points
are distributed,
and becomes large for large $R$.
As a result, we know that
the parameters $R\approx 3.2$ and $N\geq 25$ will yield
a better regularized solution.

\subsection{Irregular domain}\label{sec:cc}
As the next example, we assume the exact solution as same as
the one in the previous example
in an irregular domain
enclosed by the boundary
\[ x(\theta)=r(\theta)\cos\theta, \quad
 y(\theta)=r(\theta)\sin\theta, \qquad 0\leq\theta<2\pi \]
with the Cassini oval
\begin{equation}
 r(\theta)=
r(\theta;a,b)=a\sqrt{\cos 2\theta + \sqrt{(b/a)^4-\sin^22\theta}},
\qquad 0\leq\theta<2\pi
\label{eq:oval}
\end{equation}
in the polar coordinates
with $a=1,b=1.01$.
It is easy to show that the unit outward normal
to the boundary
is expressed as 
\[
\bm n(\theta)=\left(
\frac{r'(\theta)\sin\theta+r(\theta)\cos\theta}
{\sqrt{(r'(\theta))^2+r(\theta)^2}},
\frac{-r'(\theta)\cos\theta+r(\theta)\sin\theta }
{\sqrt{(r'(\theta))^2+r(\theta)^2}}
\right).
\]

\begin{figure}[H]
\begin{center}
\includegraphics[width=8cm]{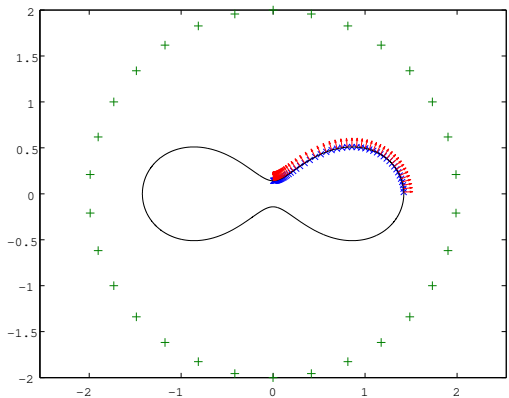}
\end{center}
\caption{Domain and distribution of points ($(R,M,N)=(2,40,10)$)}
\label{fig:ex3-dom}
\end{figure}

We distribute collocation and source points
along the boundary and the circle 
$(R\cos\theta,R\sin\theta)$ uniformly as similar as (\ref{eq:dispnt}).
Figure \ref{fig:ex3-dom} shows
the domain, the unit outward normal, the collocation and
the source points for example.

\begin{figure}[H]
\begin{center}
\includegraphics[width=8cm]{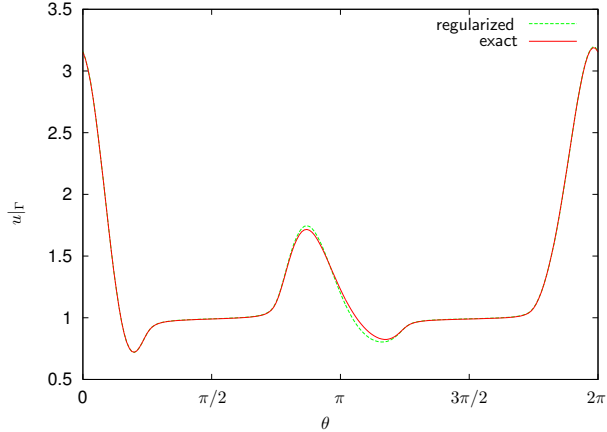}
\end{center}
\caption{The regularized solution $u_N$ versus the exact $u$ 
($10\%$ noise level; $(R,M,N)=(2,5200,30)$)}
\label{fig:ex3-2}
\end{figure}

The relative noise level of the Cauchy data
is assumed to be $10\%$ ($\delta=0.1$).
Let $(R,M,N)=(2,5200,30)$.
The optimal regularization parameter can be found as
$\alpha=2.3493\times 10^{-2}$ by using the L-curve.
Figure \ref{fig:ex3-2} shows the regularized solution
on the boundary with respect to the optimal regularization
parameter.
From the result, it is concluded that
even if the boundary of the domain is complicated
and the noise level of the Cauchy data is higher,
the regularized solution is in very good agreement with
the exact one.

\subsection{Problems with singular points}
We consider two problems whose solutions have singular points.

We first assume the exact solution
\[ u(x,y)=\log\sqrt{(x-0.2)^2+y^2}-\log\sqrt{(x+0.2)^2+y^2} \]
in the annulus domain
\[ \Omega=\{ (x,y) : 0.5^2<x^2+y^2<1 \} \]
with the outer and the inner boundaries
\[ \Gamma_{\rm out}=\{(x,y): x^2+y^2=1 \},\qquad
%and the inner boundary
 \Gamma_{\rm in}=\{(x,y): x^2+y^2=0.5^2 \}. \]
The exact solution $u$ has two singular points
at $(x,y)=(-0.2,0),(0.2,0)$.

We now assume that the exact potential $u$ is unknown.
From the Cauchy data given
on the fourth part of the outer boundary $\Gamma_{\rm out}$,
defined by
\[ \Gamma_1:=\{(x,y): x^2+y^2=1,\ x>0,\ y>0\}, \]
we identify a boundary value on the rest of the boundary
$(\Gamma_{\rm out}\setminus\Gamma_1)\cup\Gamma_{\rm in}$.

We distribute uniformly the collocation points
$\{\bm x_i\}_{i=1}^M\subset\Gamma_1$.
The source points 
$\{\bm \xi_j\}_{j=1}^{N}\subset\overline{\Omega}^c$
are uniformly distributed along
two circles whose centers are 0 and
radii are $R_{\rm out}$ and $R_{\rm in}$, respectively.

\begin{figure}[H]
\centering
\includegraphics[width=8cm]{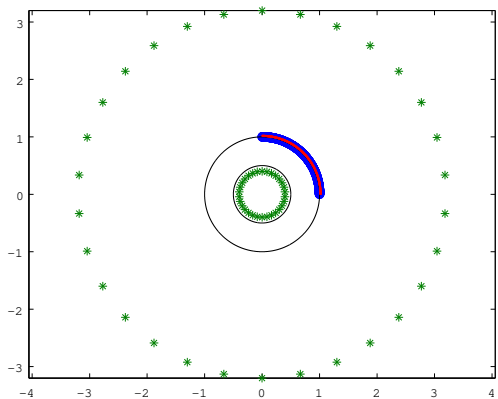}
\caption{Domain and distribution of points
 ($(R_{\rm out},R_{\rm in},M,N)=(3.2,0.4,600,60)$)}
\label{fig:ex4-dom}
\end{figure}

\begin{figure}[H]
\centering
\includegraphics[width=8cm]{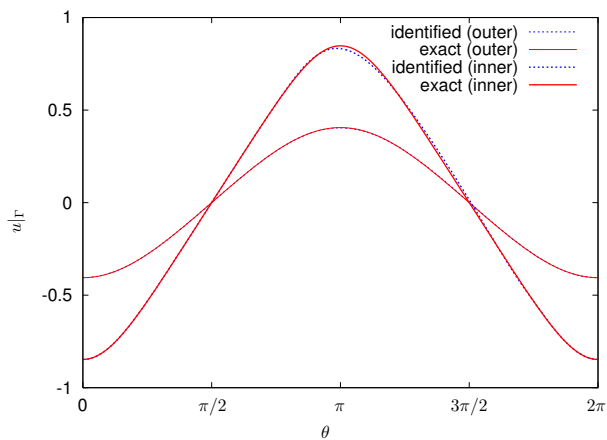}
\caption{The regularized solution $u_N$ versus the exact $u$ 
(no noise)}
\label{fig:ex4sol}
\end{figure}

In the first case, 
the exact Cauchy data is assumed to be given.
Let $(R_{\rm out},R_{\rm in},M,N)=(3.2,0.4,600,60)$ (Figure
\ref{fig:ex4-dom}).
Figure \ref{fig:ex4sol} shows the identified solution
on the outer and the inner boundaries.
We can see that the identified solution
is in very good agreement with the exact one
in spite of the solution with singular points.

\begin{figure}[H]
\centering
\includegraphics[width=8cm]{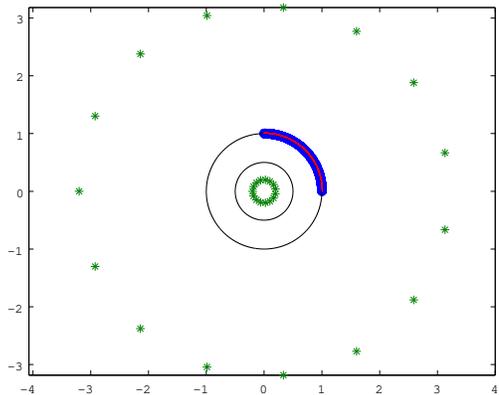}
\caption{Domain and distribution of points
 ($(R_{\rm out},R_{\rm in},M,N)=(3.2,0.2,5200,30)$)}
\label{fig:ex4-dom2}
\end{figure}

\begin{figure}[H]
\centering
\includegraphics[width=8cm]{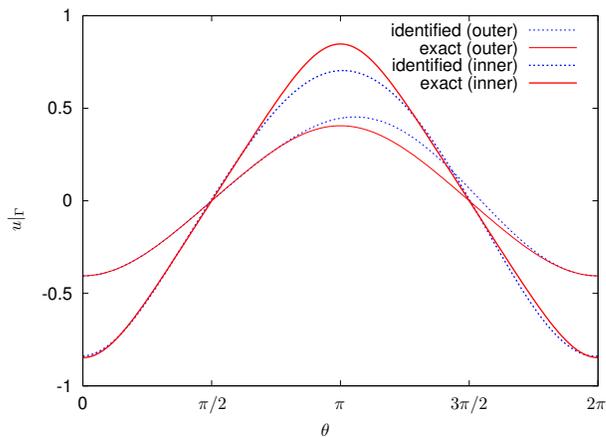}
\caption{The regularized solution $u_N$ versus the exact $u$ 
($5\%$ noise level)}
\label{fig:ex4sol-2}
\end{figure}

In the second case,
the relative noise level of the Cauchy data
is assumed to be $5\%$ ($\delta=0.05$).
Let $(R_{\rm out},R_{\rm in},M,N)=(3.2,0.2,5200,30)$
(Figure \ref{fig:ex4-dom2}).
The optimal regularization parameter can be found as
$\alpha=6.3932\times 10^{-3}$ by using the L-curve.
Figure \ref{fig:ex4sol-2} shows the regularized solution
on the outer and the inner boundaries
 with respect to the optimal regularization
parameter.
From this result, we know that the regularized solution
is acceptable.

\begin{figure}[H]
\centering
\includegraphics[width=8cm]{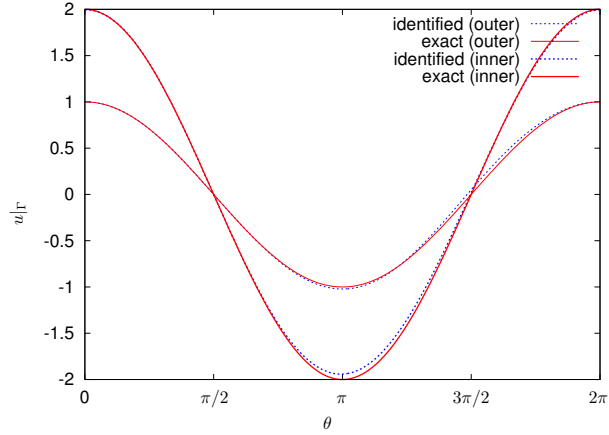}
\caption{The regularized solution $u_N$ versus the exact $u$ 
($5\%$ noise level)}
\label{fig:ex5}
\end{figure}

As another example, we assume that the exact solution is
given by
\[ u(x,y)=\frac x{x^2+y^2} \]
in the same domain as above.
The Cauchy data with $5\%$ noise level are prescribed
on the same part of the boundary as above.
Let $(R_{\rm out},R_{\rm in},M,N)=(3.2,0.05,5200,30)$.
Figure \ref{fig:ex5} shows the regularized solution
on the outer and the inner boundaries
 with respect to the optimal regularization
parameter $\alpha=1.5088\times 10^{-3}$.
The accuracy of the regularized solution is quite good.

\section{Conclusions}\label{sec:cnc}
We consider using MFS
as the numerical method for the Cauchy problem of the Laplace
equation.
Since MFS is a messless method, we can easily
treat a complicated boundary.
This paper proposes
a direct method instead of an iterative one.
The Tikhonov regularization can find a stable solution.
The L-curve automatically gives a suitable regularization parameter,
which coincides with the optimal one shown in the numerical
 experiments.
Hence, after setting the parameters $(R,M,N)$
the optimal regularized solution can be obtained
quickly and automatically.
Moreover, our numerical method can successfully
solve even a problem whose solution has singular points.

The following is the guideline for choosing better parameters $(R,M,N)$:
The collocation points should be distributed as many as possible
compared with the source points.
There is no value increasing the number of source points $N$.
It is enough for $N\approx 30$
 to obtain a better solution.
The radius of the circle where source points are distributed
should be small like $R\approx 3$,
since the stability is more important than the accuracy
in this inverse problem.

In conclusion,
the numerical method proposed in this paper is applicable
for solving a problem in a complicated domain
with the Cauchy data that contains large noises
even with a noise level of 10\%.
This method is also effective for solving even a problem
with singular points.

\section*{Acknowledgements}
The authors gratefully acknowledge the financial support
of the National Science Council of Taiwan through
the grant No.\ NSC96--2811--E--002--047.
They express their gratitude
to Professor K.~Onishi and Dr.\ K.~Shirota at
Ibaraki University, Japan,
for beneficial suggestions.

\end{document}